\documentclass[11pt, a4paper]{amsart}

\usepackage[english]{babel}
\usepackage{amsmath,enumerate, amsfonts, amssymb,amsthm, xypic}
\usepackage[dvips]{graphics} 

\newtheorem{thm}{Theorem}[section]
\newtheorem{lemma}[thm]{Lemma}
\newtheorem{prop}[thm]{Proposition}

\newenvironment{demo}{\noindent{\it Proof.}\,}{\begin{flushright}
    \,$\Box$ \smallskip \end{flushright}}

\theoremstyle{definition}
\newtheorem{defi}[thm]{Definition}

\newtheorem{ex}[thm]{Example}
\newtheorem{rmk}[thm]{Remark}
\newtheorem{Q}[thm]{Question}

\DeclareMathOperator{\ord}{ord}
\DeclareMathOperator{\jac}{jac}
\DeclareMathOperator{\mult}{mult}

\title{Zeta functions and blow-Nash equivalence}

\author{Goulwen Fichou}

\address {Institut Math\'ematiques de Rennes, Universit\'e de Rennes 1, Campus de Beaulieu, 35042 Rennes Cedex, France}

\email{goulwen.fichou@univ-rennes1.fr}
\begin{document}

\begin{abstract} We propose a refinement of the notion of
  blow-Nash equivalence between Nash function germs, which has
  been introduced in \cite{fichou} as an analog in the Nash setting of the
  blow-analytic equivalence defined by T.-C. Kuo \cite{Kuo}. The new
  definition is more natural and geometric. Moreover, this equivalence
  relation still does not admit moduli for a Nash family of isolated singularities.
But if the zeta functions
  constructed in \cite{fichou} are no longer invariants for this new
  relation, however, thanks to a Denef \& Loeser formula coming
  from motivic integration in a Nash setting, we managed to derive new invariants for
  this equivalence relation.
 \end{abstract}

\maketitle

The classification of real analytic function germs is a difficult
topic, especially in the choice of a good equivalence relation between
germs to study. Even in the particular case when the analytic function germs
are Nash, that is they are moreover semi-algebraic, the difficulty
still remains.

In \cite{fichou}, we have defined the blow-Nash equivalence between
Nash function germs, as an approximation with algebraic data
of the blow-analytic equivalence defined by T.-C. Kuo in \cite{Kuo}.
This blow-analytic equivalence has already been studied with slightly different definitions
since the original definition of T.-C. Kuo appeared (see in particular S. Koike
\& A. Parusi\'nski \cite{KP} and T. Fukui \& L. Paunescu
\cite{FP}). Nevertheless, roughly speaking, it states that two given real analytic
function germs are equivalent if they are topologicaly equivalent and
moreover, after suitable modifications, they
become analytically equivalent. 

For the case of Nash function germs, the definition of blow-Nash
equivalence runs as follows.
Let
$f,g:(\mathbb R^d,0) \longrightarrow (\mathbb R,0)$ be Nash
function germs. They are said to be blow-Nash equivalent if there
exist two Nash modifications (we refer to definition \ref{defbN} for the precise definition)
$$\pi_f:\big( M_f,\pi_f^{-1}(0) \big) \longrightarrow (\mathbb R^d,0)
\textrm{~~and~~} \pi_g:\big( M_g,\pi_g^{-1}(0)\big) \longrightarrow (\mathbb R^d,0),$$ and a
Nash isomorphism $\phi :\big( M_f,\pi_f^{-1}(0)\big) \longrightarrow
\big( M_g,\pi_g^{-1}(0)\big) $, that is $\phi$ is a real analytic
isomorphism with semi-algebraic graph,  which
induces a homeomorphism $h$ between neighbourhoods of $0$ in $\mathbb
R^d$ such that $f=g\circ h$.

For a stronger notion of blow-Nash
equivalence, we known computable invariants, which seems to be efficient tools to distinguish
blow-Nash type \cite{fichou, fichou2}. These invariants, called zeta functions (cf. section \ref{zeta}), are constructed in a
similar way to the motivic zeta functions of Denef \& Loeser, using the virtual Poincar\'e polynomial of
arc-symmetric sets as a generalized Euler characteristic (cf. section \ref{beta}).

Nevertheless, the definition of the blow-Nash equivalence given in
\cite{fichou} is strong and technical. In particular the modifications
are asked to be algebraic, which is not natural in the Nash
setting. The weaker definition of blow-Nash equivalence introduced in
this paper is more natural and geometric. It is
closer to the definition of blow-analytic equivalence considered by
S. Koike and A. Parusi\'nski in \cite{KP}. This blow-Nash equivalence
is an equivalence relation (proposition \ref{eqrel}). For such an
equivalence relation, it is a crucial fact to prove that it has a good
behaviour with respect to family of Nash function germs. In this
direction, theorem \ref{FAM} states that a family with isolated singularities does not
admit moduli. This result is
more general that the one in \cite{fichou}, whereas the present proof is just a
refinement of the former one. We mention also in section \ref{trivial} various criteria to ensure the blow-Nash triviality of a given family.

Recently, invariants for this kind of equivalence relations have been
introduced (see \cite{fukui} for a survey). In particular, we defined in
\cite{fichou} zeta functions, following ideas coming from motivic
integration \cite{DL1}, which are defined via the virtual Poincar\'e
polynomial \cite{MCP}. 

Unfortunately, if this
definition of the blow-Nash equivalence in this paper is more natural and geometric, the zeta functions are no longer
invariants in general.
However, one can derived new invariants from these zeta
functions by evaluating its coefficients, which are rational functions
in the
indeterminacy $u$ with coefficients in $\mathbb Z$ at
convenient values (cf. theorem \ref{main}). As a key ingredient, we
generalize the Denef \& Loeser formulae, that express the zeta
functions in terms of a modification, in the setting of Nash
modifications (see part \ref{proofDL}).

As a application, we manage to distinguish the
blow-Nash type of some Brieskorn polynomials whose blow-analytic type is not even known! 

\vskip 2mm

{\bf Acknowledgements.} The author wish to thank T. Fukui, S. Koike and A. Parusi\'nski for valuable discussions on the subject.

\section{Blow-Nash equivalence}\label{bN}
\subsection{}
Let us begin by stating the definition of the blow-Nash equivalence between
Nash function germs that we consider in this paper. It consists of a natural adaptation of the blow-analytic equivalence defined by T.-C. Kuo (\cite{Kuo}) to the Nash framework.
\begin{defi}\label{defbN}\begin{flushleft}\end{flushleft}
\begin{enumerate}
\item Let $f:(\mathbb R ^d,0) \longrightarrow (\mathbb R,0)$ be a Nash
  function germ. A Nash modification of $f$ is a proper
  surjective Nash map $\pi:\big( M,\pi^{-1}(0) \big)
  \longrightarrow (\mathbb R^d,0)$ whose complexification $\pi^*$ is an isomorphism except on some thin subset of $M^*$, and 
such that $f \circ \pi$ has only normal crossings.
\item \label{point2} Two given germs
of Nash functions $f,g:(\mathbb R ^d,0) \longrightarrow (\mathbb R,0)$ are said to be blow-Nash equivalent if there
exist two Nash modifications $$\sigma_f~:~\big(M_f,\sigma_f^{-1}(0)\big)  \longrightarrow
(\mathbb R ^d,0) \textrm{ and } \sigma_g:\big(M_g,\sigma_g^{-1}(0)\big)  \longrightarrow
(\mathbb R ^d,0),$$ and a Nash isomorphism $\Phi$
between semi-algebraic neighbourhoods $\big(M_f,\sigma_f^{-1}(0)\big)$ and $\big(M_g,\sigma_g^{-1}(0)\big)$ which
induces a homeomorphism $\phi :(\mathbb R ^d,0) \longrightarrow
(\mathbb R^d,0)$ such that the diagram
$$\xymatrix{\big(M_f,\sigma_f^{-1}(0)\big) \ar[rr]^{\Phi} \ar[d]_{\sigma_f}& &\big(M_g,\sigma_g^{-1}(0)\big) \ar[d]^{\sigma_g}\\
            (\mathbb R ^d,0)       \ar[rr]^{\phi} \ar[dr]_f        & &           (\mathbb R ^d,0) \ar[dl]^g \\
         &   (\mathbb R,0) & } $$
is commutative.
\end{enumerate}
\end{defi}

\begin{rmk}\begin{flushleft}\end{flushleft}
\begin{enumerate}
\item Let us specify some classical terminology (see \cite{fukui} for example). Such a homeomorphism $\phi$ is called a blow-Nash
  homeomorphism. If, as in \cite{fichou}, we ask moreover $\Phi$ to preserve the
  multiplicities of the jacobian determinant along the exceptionnal divisors
of the Nash modifications
  $\sigma_f,~\sigma_g$, then $\Phi$ is called a blow-Nash isomorphism.

Nota that there exist blow-Nash homeomorphisms which are not blow-Nash isomorphisms (see \cite{fukui}).

\item In \cite{fichou}, we consider a more particular notion of
  blow-Nash equivalence. Namely, the Nash modifications were replaced
  by proper algebraic birational morphisms and the blow-Nash homeomorphism was moreover asked to be a blow-Nash isomorphism. The definition \ref{defbN}
  is more natural since all the data are of Nash class.
\end{enumerate}
\end{rmk}

The proof of the following result is the direct analog of the
corresponding one in \cite{Kuo}.

\begin{prop}\label{eqrel} The blow-Nash equivalence is
  an equivalence relation between Nash function germs.
\end{prop}

\begin{proof} The point is the transitivity property. Let
  $f_1,f_2,f_3:(\mathbb R ^d,0) \longrightarrow (\mathbb R ,0)$ be
Nash function germs such that $f_1 \sim f_2$ and $f_2 \sim f_3$. Let
$\sigma _1, \sigma _2$ and $\sigma _2', \sigma_3'$ be Nash
modifications, and $\phi, \phi'$ be homeomorphisms like in definition \ref{defbN} for  $f_1,f_2$ and
$f_2,f_3$ respectively. The fiber product $M$ (respectively $M'$) of
$\phi \circ \sigma_1$ and $\sigma_2$ (respectively $\phi' \circ
\sigma_2'$ and $\sigma_3'$) gives suitable Nash modifications of $
(\mathbb R ^d,0)$. Taking once more the fiber product $M''$ of $M$ and $M'$
solves the problem since the compositions of the projections with the
initial modifications $\sigma_1$ and $\sigma'_3$ remain Nash
modifications for $f_1$ and $f_3$.

$$
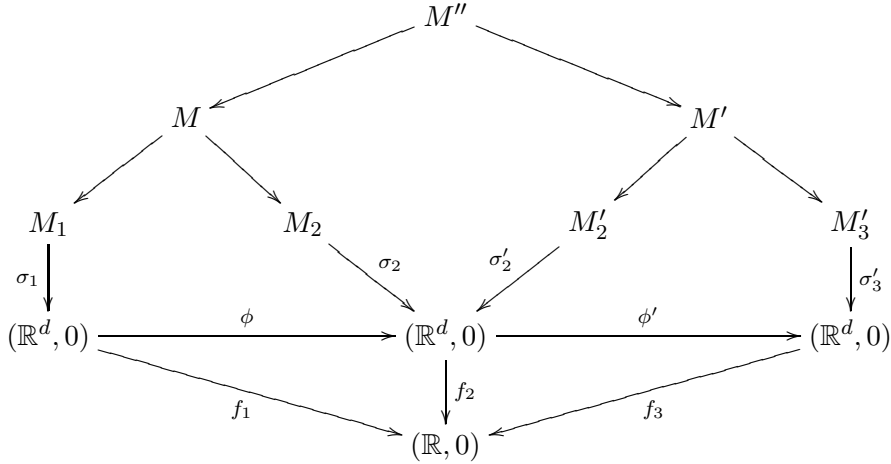
\xymatrix{        &&& M'' \ar[dll]\ar[drr] &&&       \\
& M     \ar[dl]\ar[dr]  &&& &  M'   \ar[dl]\ar[dr] &&  \\
M_1 \ar[d]_{\sigma _1} & & M_2\ar[dr]^{\sigma _2} & & M_2'\ar[dl]_{\sigma _2'} & & M_3'\ar[d]^{\sigma _3'} \\
            (\mathbb R ^d,0)       \ar[rrr]^{\phi} \ar[drrr]_{f_1}
            & &    &       (\mathbb R ^d,0) \ar[d]^{f_2}\ar[rrr]^{\phi'}  & &&
            (\mathbb R ^d,0) \ar[dlll]^{f_3}\\
    &  &   &   (\mathbb R,0) & & & } $$
\end{proof}

\begin{rmk} Note that for the blow-Nash equivalence considered in
  \cite{fichou}, we had to consider the equivalence relation generated
  by a similar relation. This difficulty came from the fact that the fiber
  product of an algebraic map and a Nash map needs not to be
  algebraic. The point here is that the fiber product of Nash maps
  still remains in the Nash class.
\end{rmk}

The question of moduli is a natural and crucial issue when one studies
an equivalence relation between germs. The following theorem
states that
the blow-Nash equivalence has a good behaviour with respect to family
of Nash function germs. More precisely, the blow-Nash equivalence does not admit
 moduli
for a Nash family of
Nash function germs with an
isolated singularity. Let's $P$ denote the cuboid $[0,1]^k$ for an integer $k$. 

\begin{thm}\label{FAM}  Let $F: (\mathbb R^d,0)\times P \longrightarrow (\mathbb R,0)$
  be a Nash map and assume that $F(.,p): (\mathbb
  R^d,0)\longrightarrow \mathbb R$ has an isolated singularity at $0$
  for each $p\in P$. 

Then the family $F(.,p)$, for $p\in P$, consists of a finite number of
  blow-Nash equivalence classes.
\end{thm}

\begin{rmk} The proof of theorem \ref{FAM} can be performed in a
  similar way to the one in
\cite{fichou}, even if the result is more general here. Indeed, we had
to restrist the study in \cite{fichou} to particular Nash families,
that is falilies which
admit, a resolution of the singularities, an algebraic modification. But,
if we allow the modifications to become Nash, the Hironaka's resolution
of singularities provides us suitable Nash modifications \cite{hiro}.
\end{rmk}


\subsection{Blow-Nash triviality}\label{trivial}

In view of classification problems, a worthwhile issue is to give criteria for a Nash family to consist of a unique blow-Nash class. In particular, one says that a Nash family $F: (\mathbb R^d,0)\times P \longrightarrow (\mathbb R,0)$ is blow-Nash trivial if there exist a Nash modification $\sigma~:~(M,E)  \longrightarrow (\mathbb R ^d,0)$, a $t$-level preserving homeomorphism $\phi: (\mathbb R ^d,0)\times P \longrightarrow (\mathbb R ^d,0)\times P$ and a $t$-level preserving Nash isomorphism $\Phi: (M,E)\times P \longrightarrow (M,E)\times P$ such that the diagram 
$$\xymatrix{    (M,E)\times P \ar[rr]^{\sigma \times \textrm{id}}\ar[d]_{\Phi} &  & (\mathbb R ^d,0)\times P \ar[rr]^{(x,p)\mapsto F(x,0)}\ar[d]^{\phi}  &&   (\mathbb R,0)\ar[d]^{\textrm{id}}\\
                (M,E)\times P \ar[rr]^{\sigma \times \textrm{id}}         && (\mathbb R ^d,0)\times P \ar[rr]^{(x,p)\mapsto F(x,p)}      &  &   (\mathbb R,0)}$$
is commutative.

Below, we mention sufficient conditions to ensure the blow-Nash triviality of a given family, that are analogs of corresponding results concerning blow-analytic equivalence (\cite{FP}, \cite{FY}). Moreover their proof (cf. remark \ref{bntrivial}) is a direct consequence of the one of theorem \ref{FAM}.

Let us introduce some terminology before stating the first result, which is inspired by the main theorem of \cite{FY}.
For an analytic function germ $f:(\mathbb R^d,0) \longrightarrow (\mathbb R,0)$, denote by $\sum_I c_I x^I$ its Taylor expansion at the origin, where $x^I=x_1^{i_1}\ldots x_d^{i_d}$, $I=(i_1,\ldots,i_d)$. The Newton polygon of $f$ is the convex hull of the union of the sets $I+\mathbb R_+^d$, for those $|I|$ such that $c_I\neq 0$. For a face $\gamma$ of this polyhedron, we put $f_{\gamma}(x)=\sum_{I\in \gamma} c_I x^I$. The germ $f$ is said to be non-degenerate, with respect to its Newton polygon, if the only singularities of $f_{\gamma}$ are concentrated in the coordinate hyperplanes, for any compact face $\gamma$ of the Newton polygon. Finally, one says that a given face is a coordinate face if it is parallel to some coordinate hyperplane.

\begin{prop}\label{tri1} Assume that the Newton polygon of $F(.,p)$ is independent of $p \in P$, non-degenerate for each $p \in P$, and moreover assume that $\big( F(.,p) \big)_{\gamma}$ in independent of $p \in P$ for any non-compact and non-coordinate face $\gamma$ of the Newton polygon. Then the family $\{F(.,p)\}_{p \in P}$ is blow-Nash trivial.
\end{prop}

The second result is inspired by the main theorem in \cite{FP}. Consider the Taylor expansion $F(x,p)=\sum_I c_I(p) x^I$ of $F$ at the origin of $\mathbb R^d$. For an $d$-uple of positive integers $w=(w_1,\ldots,w_d)$, we set $H_i^{(w)}(x,p)=\sum_{I: ~|I|_w=i} c_I x^I$, where $|I|_w=i_1w_1+\cdots+i_dw_d$. Denote by $m$ the smallest integer $i$ such that $H_i^{(w)}(x,p)$ is not identically equal to 0.

\begin{prop}\label{tri2} If there exists an $d$-uple of positive integers $w$ such that $H_m^{(w)}(x,p)$ has an isolated singularity at the origin of $\mathbb R^d$ for any $p\in P$, then the family $\{F(.,p)\}_{p \in P}$ is blow-Nash trivial.
\end{prop}

\begin{ex}(\cite{FP}) Let  $F:(\mathbb R^3,0)\times \mathbb R  \longrightarrow (\mathbb R,0)$ be the Brian\c con-Speder family, namely 
$$F(x,y,z,p)=z^5+py^6z+xy^7+x^{15}.$$
This family is weighted homogenous with weight $(1,2,3)$ and weighted degree $15$. Moreover it defines and isolated singularity at the origin for $p\neq p_0=-\frac{15^{\frac{1}{7}}\frac{7}{2}^{\frac{4}{5}}}{3}$. Therefore the Brian\c con-Speder family is blow-Nash trivial over all interval that does not contain $p_0$.
\end{ex}

\begin{rmk}\label{bntrivial} The proof of these triviality results in the blow-analytic case is based on the integration along an analytic vector field defined on the parameter space, and that can be lifted through the modification. The flow of the lifted vector field gives the trivialisation upstairs. Moreover the assumptions made enable to choose, as a modification, a toric modification that has an unique critical value at the origin of $\mathbb R ^d$. Therefore the trivialisation upstairs induces a trivialisation at the level of the parameter space.

However, by integration along a Nash vector field, one needs not keep Nash data, and therefore the same method as in the blow-analytic case does not apply in the situation of propositions \ref{tri1}, \ref{tri2}. Nevertheless, one can replace this integration by the following argument (exposed with details in \cite{fichou}).
First, resolve the singularities of the family via the relevant toric modification as in \cite{FP}, \cite{FY}. Then, trivialise the zero level of the function germs with the Nash Isotopy Lemma \cite{FKS}. Finally, trivialise the $t$-levels, $t \neq 0$, via well-choosen projections that can be proven to be of blow-Nash class.
\end{rmk}
\section{Zeta functions}

In this section, we recall the definition of
the naive zeta function of a Nash function germ (as it is defined in
\cite{fichou}). Then we prove the so-called Denef \& Loeser
formula for such a zeta function in terms of a Nash
modification. This result is new and requires to generalize the change of variables
formula in the theory of motivic integration to the Nash setting.


\subsection{Virtual Poincar\'e polynomial of arc-symmetric sets}\label{beta}
Arc-symmetric sets have been introduced by K. Kurdyka \cite{KK1} in
1988 in
order to study ``rigid components'' of real algebraic varieties. The category of arc-symmetric sets contains the real algebraic
varieties and, in some sense, this category has a better behaviour
that the one of real algebraic varieties, maybe closer to complex
algebraic varieties. For a detailed treatment of arc-symmetric sets, we refer
to \cite{fichou}. Nevertheless, let us precise the definition of such sets.

We fix a compactification of $\mathbb R^n$, for instance $\mathbb R^n
\subset \mathbb P^n$.
\begin{defi} Let $A \subset \mathbb P^n$ be a semi-algebraic set. We say
  that $A$ is arc-symmetric if,
 for every real analytic arc $\gamma :
  ]-1,1[ \longrightarrow  \mathbb P^n$ such that $\gamma
  (]-1,0[) \subset A$, there exists $\epsilon  > 0$ such that $\gamma
  (]0,\epsilon [) \subset A$.
\end{defi}

One can think about arc-symmetric sets as the
biggest category, denoted $\mathcal {AS}$, stable under boolean operations and containing the compact
real algebraic varieties and their connected components.

In particular, the following lemma specifies what the nonsingular
arc-symmetric sets are.
Note that by an isomorphism between arc-symmetric sets, we mean a birational
map containing the arc-symmetric sets in the support. Moreover, a
nonsingular arc-symmetric set is an arc-symmetric whose intersection
with the singular locus of its Zariski closure is empty.

\begin{lemma}(\cite{fichou}) Compact nonsingular arc-symmetric sets are isomorphic to
  unions of
  connected components of compact nonsingular real algebraic varieties.
\end{lemma}
A Nash isomorphism
between arc-symmetric sets $A_1,A_2 \in \mathcal{AS}$ is the
restriction of an analytic
and semi-algebraic isomorphism between compact semi-algebraic and real analytic sets $B_1,B_2$
containing $A_1,A_2$ respectively.
Generalized Euler characteristics for arc-symmetric sets are the invariants, under Nash isomorphisms, which enable to
give concrete measures in the theory of motivic integration. A generalized Euler characteristic is defined as follows.

An additive map on $\mathcal {AS}$ with values
  in an abelian group is a map $\chi$ defined on $\mathcal {AS}$
  such that
\begin{enumerate}
\item for arc-symmetric sets $A$ and $B$ which are Nash isomorphic, $\chi(A)=\chi(B)$,
\item for a closed arc-symmetric subset $B$ of $A$, 
  $\chi(A)=\chi(B)+\chi(A \setminus B)$.
\end{enumerate}
If moreover $\chi$ takes its values in a commutative ring and satisfies $\chi(A\times B)=\chi(A) \cdot \chi(B)$
for arc-symmetric sets $A,B$, then we say that $\chi$ is a generalized
Euler characteristic on $\mathcal {AS}$.

In \cite{fichou} we proved:

\begin{prop}\label{beta-nash} There exist additive maps
  on $\mathcal {AS}$ with values
  in $\mathbb Z$, denoted $\beta_i$ and called
  virtual Betti numbers, that coincide with the classical
  Betti numbers $\dim H_i(\cdot, \frac {\mathbb Z}{2 \mathbb Z})$ on
  the connected component of compact nonsingular real algebraic varieties.

Moreover $\beta(\cdot)=\sum_{i \geq 0} \beta_i(\cdot)u^i$ is
a generalized Euler characteristic on
$\mathcal {AS}$, with values in $\mathbb Z [u]$.
\end{prop}

\begin{ex} If $\mathbb P^k$ denotes the real projective space of dimension
  $k$, which is nonsingular and compact, then $\beta(\mathbb
  P^k)=1+u+\cdots+u^k$. Now, compactify the affine line $\mathbb  A_{\mathbb R}^1$ in $\mathbb P^1$
  by adding one point at the infinity. By additivity $\beta(\mathbb  A_{\mathbb R}^1)=\beta(\mathbb P^1)-\beta(point)=u,$ and so $\beta(\mathbb  A_{\mathbb R}^k)=u^k$.
\end{ex}

\begin{rmk}\label{rmk}\begin{flushleft}\end{flushleft}
\begin{enumerate}
\item The virtual Poincar\'e polynomial is not a topological invariant
  (cf \cite{MCP}).
\item The virtual Poincar\'e polynomial $\beta$ respects the dimension of arc-symmetric sets: for $A \in \mathcal
  {AS}$,
$\dim(A)=\deg \big(\beta(A)\big).$
In particular, it assures us that a
nonempty arc-symmetric set has a nonzero value under the virtual
Poincar\'e polynomial.
\item \label{chi-beta} By evaluating $u$ at $-1$, one recover the classical Euler
  characteristic with compact supports (\cite{fichou,MCP}).
\end{enumerate}
\end{rmk}

\subsection{Zeta functions}\label{zeta}
The zeta functions of a Nash function germ are defined by
taking the value, under the virtual Poincar\'e polynomial, of certain
sets of arcs related to the germ.

Denote by $\mathcal L$ the space of formal arcs at the
origin $0 \in \mathbb R ^d$, defined by:
$$\mathcal L=\mathcal L(\mathbb R ^d,0)= \{\gamma : (\mathbb R,0) \longrightarrow (\mathbb R ^d,0)
:\gamma \textrm{ formal}\},$$
and by $\mathcal L_n$, for an integer $n$, the space of arcs truncated at the order $n+1$:
$$\mathcal L_n=\{ \gamma (t)=a_1t+a_2t^2+ \cdots a_nt^n,~a_i \in \mathbb R ^d\}.$$
Let $\pi_n:\mathcal L \longrightarrow \mathcal L_n$ be
the truncation morphism.

Consider a Nash
function germ $f:(\mathbb R ^d,0) \longrightarrow (\mathbb R,0)$. We define the naive zeta function $Z_f(u,T)$ of $f$ as the
following element of $\mathbb Z[u,u^{-1}][[T]]$: 
$$Z_f(u,T)= \sum _{n \geq 1}{\beta (\mathcal X_n)u^{-nd}T^n},$$
where $\mathcal X_n$ is composed of those arcs that, composed with
$f$, give a series with order $n$:
$$\mathcal X_n =\{\gamma \in  \mathcal L_n: ord(f\circ \gamma) =n \}=\{\gamma \in  \mathcal L_n: f\circ \gamma (t)=bt^n+\cdots,
b\neq 0 \}.$$
Similarly, we define zeta functions with signs by
$$Z_f^+(u,T)= \sum _{n \geq 1}{\beta (\mathcal X_n^+)u^{-nd}T^n},~~~Z_f^-(u,T)= \sum _{n \geq 1}{\beta (\mathcal X_n^-)u^{-nd}T^n}$$
where 
$$\mathcal X_n^{\pm} =\{\gamma \in  \mathcal L_n: f\circ \gamma (t)=\pm t^n+\cdots\}.$$
Remark that $\mathcal X_n,~\mathcal X_n^{\pm}$, for $n \geq 1$, are constructible subsets of
$\mathbb R^{nd}$, hence belong to $\mathcal{AS}$.

In \cite{fichou}, we prove that these zeta functions are invariants
for the stronger notion of blow-Nash equivalence (with blow-Nash \textit{iso}morphism). Adapted to the
present case, what we will prove is:

\begin{prop}\label{inv-zeta} Let $f,g:(\mathbb R ^d,0) \longrightarrow (\mathbb R,0)$ be germs
of Nash functions. If $f$ and $g$ are blow-Nash equivalent via a
blow-Nash isomorphism, then
$$Z_f(u,T)=Z_g(u,T),~~Z_f^{\pm}(u,T)=Z_g^{\pm}(u,T).$$ 
\end{prop}

\begin{rmk}\label{rmk}\begin{flushleft}\end{flushleft}
\begin{enumerate} 
\item We do not know whether or not the zeta functions are invariant for the blow-Nash equivalence.
\item This result is a step toward the resolution of the main issue of the paper (theorem \ref{main}): which informations can we preserve, at the level of zeta functions, with only a blow-Nash \textit{homeo}morphism instead of a blow-Nash \textit{iso}morphism.

\item Note that if the modifications appearing in the definition of the
  blow-Nash equivalence of $f$ and $g$ are moreover algebraic, the
  result is precisely
  the one in \cite{fichou}. So what we have to justify here is that Nash
modifications are allowed. The key point is the Denef \& Loeser formula (cf. next section).
\end{enumerate}
\end{rmk}

\subsection{Denef \& Loeser formulae for a Nash modification}\label{proofDL}
The key ingredient of the proof of proposition \ref{inv-zeta}, and that
will be crucial in section \ref{sect3} also, is the following Denef \&
Loeser formulae which express the zeta functions of a Nash function
germ in terms of a modification of its zero locus.
First, we state the case of the naive zeta function.

\begin{prop}\label{DLform}(Denef \& Loeser formula) Let $\sigma:\big(M,\sigma^{-1}(0)\big)  \longrightarrow (\mathbb R ^d,0)$ be a
Nash modification of $\mathbb R^d$ such that $f \circ \sigma$ and the
jacobian determinant $\jac \sigma$ have only normal crossings simultaneously,
and assume moreover that $\sigma$ is an isomorphism over the
complement of the zero locus of $f$.

Let $(f \circ \sigma)^{-1}(0)= \cup_{j \in J}E_j$ be the decomposition
of $(f \circ \sigma)^{-1}(0)$
into irreducible components, and assume
that $ \sigma^{-1}(0)=\cup_{k \in K}E_k$ for some $K \subset J$.

Put $N_i=\mult _{E_i}f \circ \sigma$ and $\nu _i=1+\mult _{E_i} \jac
\sigma$, and, for $I \subset J$, denote by $E_I^0$ the set $(\cap _{i
\in I} E_i) \setminus (\cup _{j \in J \setminus I}E_j)$. Then 

$$Z_f(u,T)=\sum_{I\neq \emptyset} (u-1)^{|I|}\beta\big(E_I^0 \cap
\sigma^{-1}(0)\big) \Phi_I(T)$$

where $\Phi_I(T)=\prod_{i \in
  I}\frac{u^{-\nu_i}T^{N_i}}{1-u^{-\nu_i}T^{N_i}}.$
\end{prop}

In the case with sign, let us define first coverings of the
exceptional strata $E_I^0$ as follows.

Let $U$ be an affine open subset of $M$
such that $f \circ \sigma=u \prod_{i\in I}y_i^{N_i}$ on $U$, where $u$
is
a Nash function that does not vanish. Let us put $$R_{U}^{\pm}=\{ (x,t) \in (E_I^0 \cap U) \times \mathbb
R; t^m=\pm \frac{1}{u(x)}\},$$ where $m=gcd(N_i)$. Then the $R_{U}^{\pm}$ glue
together along the $E_I^0 \cap U$ to give $\widetilde {E_I^{0,\pm}}$.

\begin{prop}\label{DLmono} With the assumptions and notations of proposition
  \ref{DLform}, one can express the zeta functions with sign in
  terms of a Nash modification as:
$$Z_f^{\pm}(T)=\sum_{I\neq \emptyset} (u-1)^{|I|-1}\beta\big(\widetilde{E_I^{0,\pm}} \cap
\sigma^{-1}(0)\big) \prod_{i \in I}\frac{u^{-\nu_i}T^{N_i}}{1-u^{-\nu_i}T^{N_i}}.$$
\end{prop}

\begin{rmk}
The proof of propositions \ref{DLform} and \ref{DLmono} in the Nash
case run as in the algebraic one
(cf. \cite{fichou} for example, which is already an adaptation to the real case of \cite{DL1}). In particular, in the remaining of this section, we prove that we can
apply the same method. The main point is that we
dispose of a Kontsevich change of variables formula in the Nash
case. In order to prove this, the following
lemma is crucial.
\end{rmk}

\begin{lemma}\label{fibr} Let $h: \big(M,h^{-1}(0)\big)  \longrightarrow (\mathbb
  R ^d,0)$ be a proper surjective Nash map. 

Put $$\Delta_{e}=\{\gamma \in \mathcal L(M,E); \ord _t \jac h\big
(\gamma(t)\big)=e\},$$
for an integer $e\geq 1$, and $\Delta_{e,n}=\pi_n(\Delta_e)$.

For $e \geq 1$ and
  $n\geq 2n$, then $h_n(\Delta_{e,n})$ is arc-symmetric and $h_n$ is a
  piecewise trivial fibration over $\Delta_{e,n}$, where the pieces
  are arc-symmetric sets, with fiber
  $\mathbb R^e$.
\end{lemma}

As an intermediate result, note the following elementary lemma whose
proof is based
on Taylor's formula (cf. \cite{DL1}).
\begin{lemma}\label{inter} Take $e\geq 1$ and $n\geq 2e$. Then, if $\gamma_1,
  \gamma_2 \in \mathcal L (M,E)$, then if $\gamma _1 \in \Delta_e$ and
  $h(\gamma_1) \equiv h(\gamma_2)  \mod
  t^{n+1}$ then $\gamma_2 \in \Delta_e$ and $\gamma_1  \equiv \gamma_2
  \mod t^{n-e+1}$.

\end{lemma}

\begin{proof}[Proof of lemma \ref{fibr}] It follows from lemma
  \ref{inter} that $h_n$ is injective in restriction to $\Delta_{e,n}
  \cap \pi_{n-e}\big(\mathcal L(M,E)\big)$, and that $h_n\Big(\Delta_{e,n}
  \cap \pi_{n-e}\big(\mathcal
  L(M,E)\big)\Big)=h_n(\Delta_{e,n})$. Then $h_n(\Delta_{e,n})$ is
  arc-symmetric, as being the
  image by an injective Nash map of an arc-symmetric set (more precisely a constructible set).

Now, the remaining of the proof can be carried on
exactly as in \cite{DL1}.
\end{proof}

 To obtain the Kontsevich change
  of variables formula for a Nash modification, and therefore
  propositions \ref{DLform} and \ref{DLmono}, it suffices to follow the same computation as in
  \cite{fichou}. Indeed, lemma \ref{fibr} enables to apply word by word the method exposed in \cite{fichou}, just by replacing ``constructible sets'' by ``arc-symmetric sets''.

Now we can detail the proof of proposition \ref{inv-zeta}.

\begin{proof}[Proof of proposition \ref{inv-zeta}]
Let us prove the proposition in the case of the naive zeta functions.

Let $f,g:(\mathbb R^d,0) \longrightarrow (\mathbb R,0)$ be blow-Nash equivalent Nash function germs. By
definition of the blow-Nash equivalence, there exist two Nash
modifications, joined together by a commutative
diagram as in definition \ref{defbN}\ref{point2}.

By a sequence of blowings-up with smooth Nash centres, one can make
the jacobian determinants having only normal crossings. One can assume
moreover that the exceptional divisors have also only normal crossings with
the ones of the previous Nash modifications, so that we are in
situation to apply the
Denef \& Loeser formula.

Then, it is
sufficient to prove that the expressions of the zeta
functions of the germs, obtained via the Denef \& Loeser formula, coincide.
Now, the terms of the form $\beta\big(E_I^0 \cap
\sigma^{-1}(0)\big)$ are equal since the virtual Poincar\'e
polynomial $\beta$ is invariant under Nash isomorphisms
(cf. proposition \ref{beta-nash}) and the $N_i$
remain the same because of the commutativity of the diagram (cf. definition
\ref{defbN}\ref{point2}).
Finally, the $\nu_i$
coincide due to the additional assumption on the blow-Nash
homeomorphism to be a blow-Nash isomorphism.

\end{proof}

\section{Evaluating the zeta functions}\label{sect3}

In order to perform a classification of Nash function germs under blow-Nash equivalence, one needs invariants for this equivalence relation. The only ones known until now are the Fukui invariants \cite{IKK} and the zeta functions of Koike-Parusi\'nski defined with the Euler caracteristic with compact supports \cite{KP}. However, for the stronger notion of blow-Nash equivalence, the zeta functions obtained via the virtual Poincar\'e polynomial are also invariants (cf. proposition \ref{inv-zeta}).

In this section, we define new invariants for the blow-Nash
equivalence. These new invariants
are derived from the zeta functions of a Nash function
germ introduced in section \ref{zeta}. Recall that the zeta functions are formal power series in
the indeterminacy $T$ with coefficients in $\mathbb Z[u,u^{-1}]$. Then the new invariants are obtained from the zeta functions by evaluating $u$ in an appropriate way.


\subsection{Evaluate $u$ at $-1$}

To begin with, let us note that we recover the zeta functions
defined by S. Koike and A. Parusi\'nski in \cite{KP}, which has been
proven to be invariants for the blow-analytic equivalence of real
analytic function germs, by evaluating the zeta functions of section \ref{zeta}
at $u=-1$. 

Indeed, one recover the
  Euler characteristic with compact supports by evaluating the virtual
  Poincar\'e polynomial at $u=-1$ (cf. remark
  \ref{rmk}.\ref{chi-beta}). 

\begin{rmk}\label{KPsign} We recover also the zeta functions with sign in \cite{KP} of a Nash
  function germ $f$ as $-2 Z_f^{\pm}(-1,T)$. Indeed, their ones are
  defined by considering the value under the Euler characteristic with
  compact supports $\chi _c$ of the set of arcs
$$Y_n^{\pm}:=\{\gamma \in  \mathcal L_n: f\circ \gamma (t)=bt^n+\cdots,
\pm b> 0 \}.$$
But $\mathcal X_n^{\pm} \times \mathbb R_+^* \longrightarrow
Y_n^{\pm},~~(\gamma (t),a) \mapsto \gamma (at)$ is a homeomorphism,
therefore
$$\chi _c (Y_n^{\pm})=\chi _c (\mathbb R_+^* )\cdot \chi _c (\mathcal X_n^{\pm})=-2\chi _c (\mathcal X_n^{\pm}).$$
\end{rmk}

As a consequence:

\begin{prop} Let $f,g:(\mathbb R ^d,0) \longrightarrow (\mathbb R,0)$ be blow-Nash equivalent germs
of Nash functions. Then
$$Z_f(-1,T)=Z_g(-1,T),$$
and
$$Z_f^{+}(-1,T)=Z_g^{+}(-1,T),~~Z_f^{-}(-1,T)=Z_g^{-}(-1,T).$$
\end{prop}

\begin{rmk}\begin{flushleft}\end{flushleft}\label{rmkprop}
\begin{enumerate}
\item This is also a direct consequence of the proof of proposition
  \ref{inv-zeta} because by a blow-Nash homeomorphism, just the parity
of the $\nu_i$ are preserved.
\item As an application, it follows from \cite{KP} that we can state the classification of the Brieskorn polynomials of two variables $f_{p,q}=\pm x^p
 \pm y^q$, $p,q \in \mathbb N$ under blow-Nash equivalence, by using the zeta functions evaluated at $u=-1$ and the Fukui invariants. We will see another approach in section \ref{briesk2}.
\end{enumerate}
\end{rmk}

\subsection{Evaluate $u$ at $1$}
In a similar way, one can evaluate the zeta functions at $1$. In the
case of the naive zeta function, what we obtain is only zero!
Nevertheless, one can obtain finer invariants. Actually, let us decompose the naive
zeta function $Z_f(u,T)$ of a Nash function germ $f$ in the following
way:
$$Z_f(u,T)=\sum_{l \geq 1}(u-1)^l z_{f,l}(u,T),$$
where $z_{f,l}(u,T)$ is a formal power series in $T$ with coefficient
in $\mathbb Z[u,u^{-1}]$ which is not divisible by $u-1$.

Similarly, decompose the zeta functions with sign:
$$Z_f^{\pm}(u,T)=\sum_{l \geq 0}(u-1)^l z_{f,l}^{\pm}(u,T).$$
Note that here the index of the sum may begin at $0$.

By evaluating these series in $\mathbb Z[u,u^{-1}][[T]]$ at $u=1$, one
finds new invariants for the blow-Nash equivalence.

\begin{thm}\label{main} Let $f,g:(\mathbb R ^d,0) \longrightarrow
  (\mathbb R,0)$ be blow-Nash equivalent germs
of Nash functions. Then
$$z_{f,1}(1,T)=z_{g,1}(1,T),$$
$$z_{f,0}^{\pm}(1,T)=z_{g,0}^{\pm}(1,T),$$
and
$$z_{f,2}(1,T)\equiv  z_{g,2}(1,T)~ \mod 2,$$
$$z_{f,1}^{\pm}(1,T) \equiv  z_{g,1}^{\pm}(1,T)~\mod 2.$$
\end{thm}

Note that by mod 2 congruence we mean equality of the series
considered as elements in
$\frac{\mathbb Z}{2\mathbb Z}[[T]]$.
 
\begin{rmk} For $k \geq 2$, then the series $z_{f,k}^{\pm}(1,T)$ and
  $z_{f,k+1}(1,T)$ are also invariant mod 2, but unfortunately they
  just vanish!
\end{rmk}

\begin{demo} This is a consequence of the Denef \& Loeser formulae
  given in propositions \ref{DLform} and \ref{DLmono}. Let us concentrate firstly
  on the naive case.

Actually, note that
$$z_{f,1}(1,T)=\lim_{u \rightarrow 1}  \frac{Z_f(u,T)}{u-1}
\textrm{~~~and~~~} z_{g,1}(1,T)=\lim_{u \rightarrow 1}
\frac{Z_g(u,T)}{u-1},$$
that is $z_{f,1}(1,T)$ (respectively $z_{g,1}(1,T)$) is the derivative
with respect to $u$
of $Z_f(u,T)$ (respectively $Z_g(u,T)$) evaluated at $u=1$.
One can express these quotients via the Denef \& Loeser formula
(proposition \ref{DLform}). As $Z_f(u,T)$ and $Z_g(u,T)$ are divisible
by $u-1$, these quotients
coincide except the coefficients $\nu_i$, which only have the same parity. By
evaluating $u$ at $1$, we obtain the equality
 $$z_{f,1}(1,T)=z_{g,1}(1,T).$$

Similarly, $z_{f,2}(1,T)$ is the derivative of $\frac{Z_f(u,T)}{u-1}$
evaluated at $u=1$. However, the derivative of quotients of the
type $\frac{u^{-\nu}T^{N}}{1-u^{-\nu}T^{N}} $ arriving in the
expression of the Denef \& Loeser formula for $Z_f(u,T)$ are of the
form
$$-\nu \frac{u^{\nu-1}T^N}{(1-u^{-\nu}T^{N})^2}.$$
Therefore the mod 2 congruence of $z_{f,2}(1,T)$ and $z_{g,2}(1,T)$
comes from the mod 2 congruence of the different $\nu$.

One just have to repeat
the same arguments with $z_{f,0}^{\pm}(1,T)$ and $z_{f,1}^{\pm}(1,T)$ in order to complete the proof of the theorem in the cases with sign.
\end{demo}

\begin{ex} Let $f_{p,k}$ be the Brieskorn polynomial defined by $$f_{p,k}=\pm(x^p +y^{kp}
  +z^{kp}),~ p \textrm{~even},~  k\in \mathbb N.$$ 
It is not known whether two such polynomials are blow-analyticaly equivalent or not. However we prove below that for fixed $p$ and different $k$, two such polynomials are not
  blow-Nash equivalent.

Note that in \cite{fichou}, we established the analog result concerning the blow-Nash equivalence via blow-Nash \textit{iso}morphism, by using the naive zeta
functions. Actually, the naive zeta function $Z_{f_{p,k}}$ of
$f_{p,k}$ looks like
$$Z_{f_{p,k}}=(u-1) \big(
u^{-1}T^{p}+u^{-2}T^{2p}+\cdots+u^{-(k-1)}T^{(k-1)p}\big) +(u^3-1)u^{-k-2}T^{kp}$$
$$+(u-1)\big(
u^{-(k+3)}T^{(k+1)p}+u^{-(k+4)}T^{(k+2)p}+\cdots+u^{-(2k+1)}T^{(2k-1)p}\big)$$
$$+(u^3-1)u^{-2(k-2)}T^{2kp}+\cdots.$$
Now, for $p$ fixed and $k < k'$, the
$pk$-coefficient of $ Z_{f_{p,k}}$ is $(u^3-1)u^{-k-2}$ whereas the
one of $ Z_{f_{p,k'}}$ is $(u-1)u^{-k}$.
Therefore, the
$pk$-coefficient of $ z_{f_{p,k},1}$ equals $2$ whereas the
one of $ z_{f_{p,k'},1}$ is $1$, and so $f_{p,k}$ and $f_{p,k'}$ are
not blow-Nash equivalent.
\end{ex}

\subsection{Classification of two variables Brieskorn polynomials}\label{briesk2}
Effective classification of function germs under a ``blow-type''
equivalence relation is a difficult topic. In this direction, the
simplest example people tried to handle with is the one of Brieskorn
polynomials. Actually, only the
classification of two variables Brieskorn polynomials has been done completely,
under blow-analytic equivalence in \cite{KP}, and also under blow-Nash equivalence via blow-Nash \textit{iso}morphism in \cite{fichou}. In remark \ref{rmkprop}, we notice moreover that the invariants used in \cite{KP} enable to conclude also for the blow-Nash equivalence. Here we present an alternative proof using only the invariants derived from the zeta functions.

Recall that two variables Brieskorn polynomials are polynomials of the
type
$$\pm x^p \pm y^{q},~  p,q\in \mathbb N.$$
As proven in \cite{KP}, the zeta functions evaluated at
$u=-1$~(cf. remark \ref{KPsign}) enables to distinguish the
blow-Nash type except in the particular case of
$$f_k(x,y)=\pm (x^k+y^k),~~ k ~\textrm{~even}.$$

In that case, by Denef \& Loeser formulae we obtain
$$Z_{f_k}(T)=(u^2-1)\frac{T^k}{u^{2}-T^k},$$
and if $f_k(x,y)=x^k+y^k$,
$$Z_{f_k}^+(T)=(1+u)\frac{T^k}{u^{2}-T^k},~~Z_{f_k}^-(T)=0,$$
and the converse if $f_k(x,y)=+(x^k+y^k)$.

Therefore
$$z_{f_k,1}=2\frac{T^k}{1+T^k}$$
and thus $z_{f_k,1} \neq z_{f_{k'},1}$ whenever $k \neq k'$,
whereas if $k=k'$ but the signs are different, the cancellation of
$z_{f_k,1}^+$ or $z_{f_k,1}^-$ enables to distinguish $f_k$ and
$f_{k'}$.

%

As a consequence, we have proved that we can draw the classification
under blow-Nash equivalence of the Brieskorn polynomials of two
variables, by using the invariants derived from the zeta functions by
evaluation of the indeterminacy $u$. Moreover, this classification
coincides with the ones established in \cite{KP} and \cite{fichou}, that is the blow-analytic, blow-Nash via blow-Nash \textit{iso}morphism and blow-Nash type of the Brieskorn polynomials of two variables are the same.


\section{Questions}

As we have already noticed, the invariants known for the
blow-analytic equivalence (the Fukui invariants \cite{IKK}, the zeta functions of
S. Koike and A. Parusi\'nski \cite{KP}) are invariants for the
blow-Nash equivalence. However:

\begin{Q} Do the zeta functions $Z_{f}(u,T)$ of a real analytic
  function germ be invariants for the blow-analytic equivalence? Or,
  as a weaker version, do
  the invariants obtained after evaluation at $1$ be invariants for the blow-analytic equivalence?
\end{Q}

More generally, the differences between the blow-Nash equivalence and
the blow-analytic one are not known in the case of Nash function
germs or even of polynomial germs. As an example, we haved proved that the blow-analytic and the blow-Nash types of the
Brieskorn polynomials of two variables coincide. But in general:

\begin{Q} Do the blow-Nash equivalence and the blow-Nash equivalence via blow-Nash \textit{iso}morphism coincide?
\end{Q}

\begin{Q} Do the blow-Nash equivalence(s) and the
blow-analytic equivalence coincide on polynomial germs? On Nash
function gems?
\end{Q}

\end{document}